\documentstyle[11pt]{article}
\newtheorem{thm}{Theorem}
\newtheorem{defn}[thm]{Definition}
\newtheorem{lemma}[thm]{Lemma}
\newtheorem{prop}[thm]{Proposition}
\newtheorem{cor}[thm]{Corollary}

\newtheorem{conj}[thm]{Conjecture}

\newtheorem{question}[thm]{Question}

\newtheorem{obs}[thm]{Observation}

\newcommand{\beq}[1]{\begin{equation}\label{#1}}
\newcommand{\enq}[0]{\end{equation}}

\newcommand{\qed}[0]{\begin{flushright} \rule{2mm}{3mm} \end{flushright}}

\newcommand{\C}[2]{{{#1}\choose{{#2}}}}
\newcommand{\ga}[0]{\alpha }
\newcommand{\gb}[0]{\beta }
\newcommand{\gc}[0]{\gamma }
\newcommand{\gd}[0]{\delta }

\newcommand{\gL}[0]{\Lambda}

\newcommand{\gO}[0]{\Omega}

\newcommand{\gs}[0]{\sigma}

\newcommand{\0}[0]{\emptyset}

\newcommand{\ra}[0]{\rightarrow}
\newcommand{\Ra}[0]{\Rightarrow}

\newcommand{\Nn}[0]{{\bf N}}

\newcommand{\xx}[0]{{\bf x}}

\newcommand{\E}[0]{{\sf E}}

\newcommand{\A}[0]{{\cal A}}
\newcommand{\B}[0]{{\cal B}}
\newcommand{\cee}[0]{{\cal C}}

\newcommand{\h}[0]{{\cal H}}
\newcommand{\I}[0]{{\cal I}}

\newcommand{\m}[0]{{\cal M}}
\newcommand{\N}[0]{{\cal N}}
\newcommand{\Oh}[0]{{\cal O}}

\newcommand{\bn}[0]{\bigskip\noindent}
\newcommand{\mn}[0]{\medskip\noindent}
\newcommand{\nin}[0]{\noindent}

\newcommand{\sub}[0]{\subseteq}
\newcommand{\sm}[0]{\setminus}
\renewcommand{\dots}[0]{,\ldots,}

\newcommand{\uone}[0]{\underline{1}}
\newcommand{\1}[0]{{\bf 1}}

\newcommand{\da}[0]{\downarrow}

\newcommand{\sugg}[1]{}

\date{}

\begin{document}

\renewcommand{\thefootnote}{\fnsymbol{footnote}}
\title{Conditional negative association for competing urns\footnotemark }

\author{
J. Kahn and M. Neiman\\
{Rutgers University and UCLA} \\
{\footnotesize email: jkahn@math.rutgers.edu; neiman@math.ucla.edu}
}
\date{}

\footnotetext{ * Supported by NSF grant DMS0701175.}

\maketitle

\begin{abstract}
We prove conditional negative association for random variables
$\xx_j=\1_{\{|\gs^{-1}(j)|\geq t_j\}}$ ($j\in [n]:=\{1\dots n\}$),
where $\gs(1)\dots \gs(m)$ are i.i.d. from $[n]$.
(The $\gs(i)$'s are thought of as the locations of balls
dropped independently into urns $1\dots n$ according to
some common distribution, so that, for some threshold $t_j$,
$\xx_j$ is the indicator of the event that at least $t_j$
balls land in urn $j$.)
We mostly deal with the more general situation in which
the $\gs(i)$'s need not be identically distributed, proving
results which imply conditional negative association in the
i.i.d. case.
Some of the results---particularly Lemma \ref{Glemma} on graph
orientations---are thought to be of independent interest.

We also give a counterexample to a negative correlation conjecture
of D. Welsh, a strong version of a (still open) conjecture of
G. Farr.
\end{abstract}

\section{Introduction}\label{Intro}

{\em Competing urns}
refers to the experiment in which $m$ balls are dropped,
randomly and independently, into urns $1\dots n$.
Formally, we have a random $\gs:[m]\ra [n]$
(where $[m]=\{1\dots m\}$)
with the $\gs(i)$'s independent.
We then take $\xx_j$ to be the indicator for occupation of urn $j$
and are interested in the law, $\mu$, of
$(\xx_1\dots \xx_n)$ (a measure on  $\{0,1\}^n$).
In the traditional case where the balls are identical
(i.e. the $\gs(i)$'s are i.i.d.) we
call $\mu$
an {\em urn measure},
or, for emphasis, an {\em ordinary} urn measure.
More generally, setting $B_j=|\gs^{-1}(j)|$,
we may consider thresholds $t_1\dots t_n$,
and let $\xx_j $ be the indicator of
$\{B_j\geq t_j\}$;
for i.i.d. balls, we then call the law of $(\xx_1\dots \xx_n)$ a
{\em threshold} urn measure.
When the balls are not required to be identical we speak of
{\em generalized urn measures} and
{\em generalized threshold urn measures}.

We are interested in correlation properties of these measures,
but before proceeding further need to briefly recall
a few definitions.
A fuller version of the following discussion is given in \cite{KN1},
and further background and motivation may be found e.g. in \cite{Pem}.

Recall that events $\A,\B$ in a probability space
are \emph{positively
correlated}---we write $\A\uparrow \B$---if $\Pr(\A\B)\geq\Pr(\A)\Pr(\B)$,
and {\em negatively} correlated ($\A\downarrow \B$)
if the reverse inequality holds.

We will be interested in measures on finite product spaces
$\Omega = \prod_{i=1}^n\gO_i$ with each $\gO_i$ a chain
(totally ordered set), often simply $\{0,1\}$.
We use
$\m(\gO)$, or simply $\m$, for the set of probability measures
on $\Omega$, and $\m_S$ for $\m(\{0,1\}^S)$.
We will occasionally identify $\{0,1\}^S$ with $2^S$
($\{\mbox{subsets of $S$}\}$ ordered by inclusion)
in the usual way.

Recall that an event $\A\sub \Omega$
is {\em increasing} (really, nondecreasing)
if $x\geq y\in\A~\Ra x\in\A$ (where $\Omega$ is endowed with the product order),
and similarly for {\em decreasing}.
For real-valued random variables $X,Y$, write $X\da Y$ if
\beq{bna}
\mbox{$\{X \geq s\}\da\{Y \geq t\} ~~~\forall s,t \in \Re$,}
\enq
or, equivalently, if
\beq{bna'}
\mbox{$\E f(X)g(Y) \leq \E f(X) \E g(Y)$ for all increasing
$f,g:\Re \to \Re$.}
\enq
(N.B. this differs from the usage in \cite{Pem}.
Of course $X\uparrow Y$ means the reverse inequalities hold,
but we don't need this.)

Say $i\in [n]$ {\em affects} $\A\sub \Omega$ if there are $\eta\in \A$
and $\tau\in\Omega\sm\A$ with $\eta_j=\tau_j ~\forall j\neq i$,
and write $\A\perp \B$ if no coordinate affects both $\A$ and $\B$.
Then $\mu\in \m$ is {\em negatively associated}
(or has {\em negative association}; we use ``NA" for either)
if
$\A\downarrow \B$
whenever $\A,\B$ are increasing and $\A\perp \B$.
We say $\mu$ has {\em negative correlations} (or {\em is NC}) if
$\eta_i\downarrow \eta_j$ (where $\eta$ is the random string)
whenever $i\neq j$.

We are primarily concerned with {\em conditional} negative association:
$\mu\in \m$ is {\em conditionally negatively associated} (CNA) if
any measure obtained from $\mu$ by conditioning on the
values of some of the variables is NA.
(Throughout the paper we assume that any conditioning
event we consider has positive probability.)
{\em Conditional negative correlation}
(CNC) for $\mu$ is defined similarly.

When $\gO=\{0,1\}^n$,
stronger properties are obtained by demanding NC (resp. NA) for
every measure $W\circ \mu\in\m$ of the form
$$W\circ \mu(\eta) ~\propto~ \mu(\eta)\prod W_i^{\eta_i}$$
with $W=(W_1\dots W_n)\in\Re_+^n$.
(Borrowing Ising terminology, one says that $W\circ\mu$
is obtained from $\mu$ by {\em imposing an external field}.)
Then $\mu$ is said to be {\em Rayleigh} or {\em NC+}
(resp. {\em NA+}), the reference in the former case being to
Rayleigh's monotonicity law for electric
networks (see e.g. \cite{DS} or \cite{choewagner03}).

\medskip
The competing urns model was
explored in some detail by
Dubhashi and Ranjan \cite{DR}\footnote{They say ``bins"
rather than ``urns."}, who
proved {\em inter alia} that threshold urn measures are NA.
Another proof of this is given in \cite{Pem}.
Actually the argument of \cite{DR}, which proves
the stronger statement that
the (law of the) r.v.'s
\beq{xiij}
\xi_{ij}={\bf 1}_{\{\gs(i)=j\}}
\enq
is NA, does not require identical balls.
(The argument of \cite{Pem} does not work for nonidentical balls.)
The main purpose of the present note is to prove
\begin{thm}\label{Uthm}
Threshold urn measures are CNA.
\end{thm}
In contrast, as observed in \cite{KN1}, even ordinary
urn measures need not be Rayleigh
(but see the remark on $R^+$ in Section \ref{Remks}).
We don't know whether Theorem \ref{Uthm} extends to
nonidentical balls (again see Section \ref{Remks}).

\medskip
Let us quickly say
what Theorem \ref{Uthm} has to do with \cite{KN1}.
Following \cite{Pem}, we say that $\mu\in \m(\{0,1\}^n)$ is
{\em ultra-log-concave}
(ULC)
if its {\em rank sequence}, $\{r_i:=\mu(|\eta|=i)\}_{i=0}^n$
(where
$|\eta|=\sum\eta_i$), has no internal zeros and the sequence
$\{r_i/\C{n}{i}\}_{i=0}^n$ is log-concave.
A set of four conjectures from \cite{Pem} (see his Conjecture 4)
states that each of CNC, CNA, NC+ and NA+ implies ULC;
but, as shown in \cite{BBL} and \cite{KN1},
even the weakest of these (NA+ $\Ra$ ULC) is false.
Theorem \ref{Uthm} provides a more natural counterexample to
the stronger ``CNA $\Ra $ ULC," since, as observed in \cite{KN1}
(disproving another conjecture from \cite{Pem}),
urn measures need not be ULC.

\medskip
The proof of Theorem \ref{Uthm} gives something a little
more general, as follows.
Suppose that for each $j\in [n]$ we are given a sequence
$0=a_0(j)<\cdots < a_{k_j}(j)=m+1$,
and for $\gs:[m]\ra [n]$ set
\beq{xx}
\xx_j(\gs) = t ~~\mbox{iff} ~ ~a_t(j)\leq B_j
< a_{t+1}(j).
\enq
\begin{thm}\label{Uthm'}
If the $\gs(i)$'s are i.i.d. then the $\xx_j$'s in {\rm (\ref{xx})}
are CNA.
\end{thm}
Call the law of $(\xx_1\dots \xx_n)$ as in (\ref{xx})
a (generalized) {\em interval urn measure}.

\medskip
The paper is organized as follows.
Section \ref{PT1} reduces Theorem \ref{Uthm'} to either
of our two main inequalities, (\ref{mukl'}) and (\ref{mukl''}).
Each of these is valid at the level of generalized interval urns;
they are equivalent in the case of ordinary urns but not obviously so
in general
(though the argument in Section \ref{FP}
uses some interplay between the two).
It is only in the derivation of
Theorem \ref{Uthm'} from (\ref{mukl'})
that we need the $\gs(i)$'s to be i.i.d.

We give two quite different ways of getting at
these main inequalities.
Theorem \ref{MainThm} in
Section \ref{FP}  essentially restates
(\ref{mukl'}) and (\ref{mukl''})
in induction-friendly form;
the proof of the theorem given in this section is
inspired by \cite{DR}.
Section \ref{GA} takes a different approach,
based on a graph-theoretic observation,
Lemma \ref{Glemma}, that is thought to  be of
independent interest.
The lemma is used to: reprove
(\ref{mukl'}); in combination with a result from
\cite{CAPP} (Theorem \ref{CAPPULC} below), to prove a stronger,
ultra-log-concavity version of (\ref{mukl''});
and to prove
``log-submodularity"
for some classes of measures.

Finally, Section \ref{Remks} contains some discussion of the
question of whether Theorem \ref{Uthm} extends to nonidentical
balls, mentions a conjecture of G. Farr and a stronger one of D. Welsh,
and sketches a counterexample to the latter.

\mn
{\bf Some notation.}
%
For a nonnegative vector
\beq{gamma}
\gc = (\gc_{ij}:i\in [m],j\in [n]),
\enq
$A\sub [m]$ and $K\sub [n]$, the probability measure
on $K^A$
(functions from $A$ to $K$)
{\em corresponding to} $\gc$ is that given by
\beq{tauW}
\Pr(\gs) \propto W(\gs):= \prod_{i\in A}\gc_{i,\gs(i)}.
\enq
Thus the r.v.'s $\gs(i)$ are independent;
they are i.i.d. if $\gc_{ij}$ does not depend on $i$,
in which case we write simply $\gc_j$.
We also use
$\Pr^{L}$ ($L\sub [m]$) for the  measure on
$[n]^L$
corresponding to $\gc$ (so $\Pr = \Pr^{[m]}$).

\section{Setting up}\label{PT1}

Let the law of $\gs\in [n]^{[m]}$
be
given by (\ref{tauW}) and let
$\xx_1\dots \xx_n$ be as in (\ref{xx}).
Let $I \cup J \cup K$ be a partition of $[n]$ and
$t_j\in \{0\dots k_j-1\}$ for $j\in K$, and set
\beq{Q}
Q:=\{ \xx(\gs) \equiv t \textrm{ on } K\}
~~~~~~(= \{a_{t_j}(j)\leq |\gs^{-1}(j)| < a_{t_j+1}(j)~\forall j \in K\}),
\enq
$X=|\gs^{-1}(I)|$ and $Y=|\gs^{-1}(J)|$.
The main point for the proof of Theorem \ref{Uthm'} is
\beq{MP}
\mbox{$ X\da Y~$ given $Q$,}
\enq
given which we finish easily:

\mn
{\em Proof of Theorem} \ref{Uthm'}.
With notation as above, let $\A,\B\sub \gO$ be
increasing events determined by $I$ and $J$
(more precisely, by the values of the variables $\xx_j$
indexed by $I$ and $J$)
respectively.
For Theorem \ref{Uthm'} we should show $ \A\da \B $ given $Q$.
Define $f,g:\Nn \to \Re$
(where $\Nn=\{0,1,\ldots\}$)
by $f(k) = \Pr( \A | X = k)$,
$g(l) = \Pr(\B | Y = l)$.
A standard coupling argument shows that
$f$ and $g$ are increasing, whence,
according to (\ref{MP}),
$$\Pr( \A \cap \B | Q) =
\E [f(X)g(Y)|Q] \leq \E [f(X)|Q] \E [g(Y)|Q]
= \Pr(\A | Q) \Pr(\B | Q)$$
(where the first equality follows from conditional
independence of $\A$ and $\B$ given $(X,Y)$).\qed

We continue to condition on $Q$
and write $\mu_k$ for the law of $Y$ given $\{X=k\}$; that is,
\beq{mu}
\mu_k(l)=\Pr(Y=l|X=k).
\enq
We will actually
prove
\beq{mukl'}
\frac{\mu_{k+1}(l+1)}{\mu_{k+1}(l)}
\leq \frac{\mu_k(l+1)}{\mu_k(l)}
\enq
(whenever neither side is $0/0$, where we
agree that $x/0=\infty$ for $x>0$),
which is a strengthening of
(\ref{MP}) once we rule out some pathologies.
We recall the standard
\begin{defn}\label{defcvx}
$~\cee\sub \Nn^n$ is {\em convex} if
$a,c\in \cee$ and $a\leq b\leq c$ imply $b\in \cee$.
\end{defn}
It will follow from Proposition \ref{propCVX} below
that
\beq{cvx}
\mbox{${\rm supp}(\Pr) := \{(k,l):\Pr(X=k,Y=l)>0\}~$ {\em is convex.}}
\enq
%
Given this,
(\ref{mukl'})
implies that $Y$ is stochastically
decreasing in $X$---that is,
$$
\mu_{k+1}(Y\geq t)\leq \mu_k(Y\geq t) ~~
\forall k,t
$$
(the easy implication is essentially Proposition 1.2 of \cite{Pem})---which
in turn easily implies $X\downarrow Y$.

\medskip
Let
\beq{Z}
Z=|\gs^{-1}(I\cup J)|.
\enq
When the $\gs(i)$'s are i.i.d., an alternate
way to specify $X$ and $Y$ is:
let $Z$ be as in (\ref{Z}),
$X\sim \textrm{Bin}(Z,\ga)$ and $Y=Z-X$, where
$\ga = \gc_I/\gc_{I\cup J}$ (with $\gc_I=\sum_{i\in I}\gc_i$)
and
$\textrm{Bin}(Z,\ga)$ is the binomial distribution with
parameters $Z$ and $\ga$.
%

In general, for $\nu$ the law of an $\Nn$-valued r.v. $Z$ and
$\ga\in [0,1]$, let
$X=X_{\nu,\ga}\sim \textrm{Bin}(Z,\ga)$, $Y=Y_{\nu,\ga}=Z-X$
and,
for lack of a better name, say
$\nu$ is {\em binomially negatively associated} (BNA)
if
$X\da Y$ for every $\ga$.
Call a nonnegative sequence $a=(a_i)_{i=0}^{\infty}$
{\em strongly log-concave} (SLC) if 
\beq{slc}
i a_i^2 \geq (i+1) a_{i-1} a_{i+1} ~~\forall i \geq 1
\enq
(that is, $(i!a_i)_{i=0}^{\infty}$ is log-concave),
and say $\nu\in \m(\Nn)$
is SLC if the sequence
$(\nu(i))_{i=0}^{\infty}$ is.
A straightforward calculation shows that
this is equivalent to
saying that (\ref{mukl'}) holds
for any $\ga$,
$X=X_{\nu,\ga}$ and $Y=Y_{\nu,\ga}$
(and $\mu_k$ as in (\ref{mu})):
since
$$\mu_k(l) = \frac{\nu(k+l)\Pr(X=k|Z=k+l)}{\Pr(X=k)}=
\frac{\nu(k+l)\C{k+l}{k}\ga^k(1-\ga)^l}{\Pr(X=k)},$$
we may rewrite (\ref{mukl'}) as
$$
\mbox{$\nu(k+l+2)\C{k+l+2}{k+1}
\nu(k+l)\C{k+l}{k}
\leq
\nu(k+l+1)\C{k+l+1}{k+1}
\nu(k+l+1)\C{k+l+1}{k},$}
$$
which is SLC for $\nu$.
(If $\nu$ is Poisson---that is, if (\ref{slc}) holds with
equality---then $X$ and $Y$
are independent Poisson r.v.'s and the inequalities
(\ref{bna}) are equalities.)
Thus, in the i.i.d. case, (\ref{mukl'}) is equivalent to
saying that $Z$ as in (\ref{Z}) is SLC.
The latter again turns out to be true at the level of {\em generalized} urns;
that is, for any $\gc$ as in (\ref{gamma}),
$\gs\in [n]^{[m]}$ with law given by (\ref{tauW}),
$Q$ as in (\ref{Q}) and $Z$ as in (\ref{Z}),
%
\beq{mukl''}
\mbox{{\em the law of $Z$ is SLC.}}
\enq
It's also easy to see that
absence of internal zeros in $(\nu(i))$ is
equivalent to (\ref{cvx})
for $X=X_{\nu,\ga}$, $Y=Y_{\nu,\ga}$
(which, again, is given by Proposition \ref{propCVX}),
so that (\ref{mukl''}) again implies (\ref{MP})
(and Theorem \ref{Uthm'}).
It seems interesting that both (\ref{mukl'}) and
(\ref{mukl''}) are valid for generalized urns, though
the equivalence that holds for i.i.d. balls
disappears in the more general setting.

As mentioned earlier, in Section \ref{GA} we will
combine Lemma \ref{Glemma} with
a result from \cite{CAPP} to obtain an improvement of
(\ref{mukl''}):
\beq{mukl'''}
\mbox{{\em the law of $Z$ is ultra-log-concave.}}
\enq
(Recall---see following
Theorem \ref{Uthm}---this means
that the
sequence $\{\Pr(|Z|=i)/\C{m}{i}\}_{i=0}^m$ is log-concave
without internal zeros.)

\section{First proof}\label{FP}

Let $\Pr$ be the measure on $[n]^{[m]}$ corresponding
to some $\gc$ (see the end of Section \ref{Intro}), and
for $a,b\in \Nn^{n-1}$ and $k \in \Nn$,
set
\beq{psab}
p(k,a,b)=\Pr(B_n=k |B_j\in [a_j,b_j] ~\forall j\in [n-1])
\enq
(recalling that $B_j=|\gs^{-1}(j)|$).
%
\begin{thm} \label{MainThm} With notation as above,

\mn
{\rm (a)} $\displaystyle \frac{p(k+1,a,b)}{p(k,a,b)}$ is
nonincreasing in $(a,b)$, and

\mn
{\rm (b)} $\displaystyle \frac{p(k+1,a,b)}{p(k,a,b)}
\leq \frac{k}{k+1} \cdot \frac{p(k,a,b)}{p(k-1,a,b)},$

\mn
where we say nothing about the case $0/0$ and
agree that $x/0=\infty$ when $x>0$.
\end{thm}
See also Theorem \ref{DR26} in Section \ref{Remks} for a related result.

As noted earlier, part (b) of Theorem \ref{MainThm} is just
a reformulation of (\ref{mukl''}), while (a) is a mild
generalization of (\ref{mukl'})
(which in fact---see (\ref{I1})---quickly reduces to (\ref{mukl'})).
To see this, note that in (\ref{mukl'}) we may assume that each of
$I,J$ is a singleton, say $K=[n-2]$, $I=\{n-1\}$, $J=\{n\}$---formally
we could pass to
$$\gc_{ij}'=\left\{\begin{array}{ll}
\gc_{ij} &\mbox{if $j\in K=[n-2]$}\\
\sum \{\gc_{ij}:j\in I\}&\mbox{if $j=n-1$}\\
\sum \{\gc_{ij}:j\in J\}&\mbox{if $j=n$}
\end{array}
\right.$$
---and similarly in (\ref{mukl''}) we may assume $K=[n-1]$, $I=\{n\}$ and
$J=\0$.
Then Theorem \ref{MainThm}(b), which
may also be stated
$$
\mbox{\em for fixed $a,b \in \Nn^{[n-1]}$ the sequence $\{p(k,a,b)\}$
is SLC},
$$
is (up to some name changes) the same as (\ref{mukl''}),
while (\ref{mukl'}) is equivalent to
$$p(k+1,a,b)/p(k,a,b) \geq p(k+1,a',b')/p(k,a',b'),$$
where
$a_{n-1}=b_{n-1}=t$, $a'_{n-1}=b'_{n-1}=t+1$ (for some $t$)
and $a_j=a'_j$, $b_j=b'_j$ for $j\in [n-2]$.

On the other hand, the inductive proof of Theorem \ref{MainThm} employs
both the more general form of (a) and some interplay
between the two parts.

\medskip
Before proving the theorem  we note one further
consequence
and
give the promised Proposition \ref{propCVX}.
For $f ,a\in \Nn^n$,
let
$\m_f(a)=\{\gs\in [n]^{[m]}:|\gs^{-1}(j)|\in [a_j,a_j+f_j]~\forall j\in [n]\}$
and $M_f(a)=\Pr(\m_f(a))$.
Though we won't use the next result
(but see the remark following Corollary \ref{gphcor}),
it seems natural and worth mentioning.
\begin{cor}\label{corNLC}
For each $f\in \Nn^n$, $M=M_f$ satisfies the negative lattice condition:
\beq{CI1}
M(a)M(c) \geq M(a \vee c)M(a \wedge c)~~\forall a,c\in \Nn^n .
\enq
\end{cor}
This
is more or less immediate
from Theorem \ref{MainThm}
once we have the next little observation, which, as noted earlier,
also gives (\ref{cvx}) and
absence of internal zeros in the law of $Z$ in (\ref{Z}).
\begin{prop}\label{propCVX}
For any $f$ and $M_f$ as above, the support of $M=M_f$ is convex.
\end{prop}
{\em Proof.}
This will follow easily from

\mn
{\em Claim.}
For any $\gs,\tau\in [n]^{[m]}$ with $\Pr(\gs),\Pr(\tau)>0$
and $i\in [n]$ with
$|\gs^{-1}(i)| > |\tau^{-1}(i)|$, there
are $j\in [n]$ and $\rho\in [n]^{[m]}$ with $\Pr(\rho)>0$,
$|\gs^{-1}(j)| < |\tau^{-1}(j)|$ and
$$
|\rho^{-1}(k)| = \left\{\begin{array}{ll}
|\gs^{-1}(i)|-1&\mbox{if $k=i$}\\
|\gs^{-1}(j)|+1&\mbox{if $k=j$}\\
|\gs^{-1}(k)|&\mbox{if $k\in [n]\sm\{i,j\}$.}
\end{array}\right.
$$
This is a standard type of
graph-theoretic observation:
regarding
$\gs$ and $\tau$ as edge sets of
bipartite graphs on $[m]\cup [n]$ in the natural
way,\footnote{We pretend $[m]\cap [n]=\0$.}
we need 
a path with edges alternately from $\gs\sm \tau$ and $\tau\sm \gs$
that begins with a $\gs$-edge at $i$ and
ends with a $\tau$-edge at some $j$ as above.
(We then get $\rho$ by switching $\gs$ and $\tau$ on
this path.)
We omit the routine proof that such a path must exist.

To prove Proposition \ref{propCVX},
we should show that for all
distinct $a,b,c \in \Nn^n$ with $a \leq b \leq c$
and $a,c \in {\rm supp}(M)$, we also have $b \in {\rm supp}(M)$.
Of course it suffices to show this when there is some $i \in [n]$
with $b_i=c_i-1$ and $b_k=c_k \textrm{ for all } k \neq i$.
Choose $\tau\in \m(a):=\m_f(a)$ and $\gs\in \m(c)$ with
$\Pr(\tau),\Pr(\gs)>0$.
We assume $|\gs^{-1}(i)|=c_i+f_i$,
since otherwise $\gs \in \m(b)$ and we are finished.
Letting $j, \rho$ be as in the claim (note $|\tau^{-1}(i)| < c_i+f_i$),
we have $$|\rho^{-1}(i)|=b_i+f_i$$
and
$$|\rho^{-1}(j)|=|\gs^{-1}(j)|+1 \in [c_j+1, a_j+f_j] \sub [b_j,b_j+f_j],$$
whence $\rho \in \m(b)$ and $b \in {\rm supp}(M)$. \qed

\mn {\em Proof of Corollary} \ref{corNLC}.
It is easy to see (and standard) that convexity of $M$
(given by
Proposition \ref{propCVX}) implies that
it's enough to prove (\ref{CI1}) when
there are indices $i$ and $j$ with
$a_i = c_i-1$, $a_j = c_j+1$,
and $a_k=c_k$ for all
$k\neq i,j$.
In this case---assuming, w.l.o.g., that
$i=n-1$ and $j=n$---we set
$$p_1(k)=\Pr(B_n=k|B_l \in [a_l, a_l + f_l] ~\forall l\in [n-1])$$
and
$$p_2(k)=\Pr(B_n=k|B_l \in [c_l, c_l + f_l] ~\forall l\in [n-1]).$$
Then (\ref{CI1}) is
$$
\left(\sum_{k=c_n+1}^{c_n+f_n+1} p_1(k)\right)
\left(\sum_{k=c_n}^{c_n+f_n} p_2(k)\right)
\geq
\left(\sum_{k=c_n+1}^{c_n+f_n+1} p_2(k)\right)
\left(\sum_{k=c_n}^{c_n+f_n} p_1(k)\right)
$$
and follows immediately from
$$p_1(k)p_2(l) \geq p_1(l)p_2(k) ~ \textrm{{\em whenever }}k \geq l,$$
which is a consequence of
Theorem \ref{MainThm}(a)
(and Proposition \ref{propCVX}).
\qed

We now assume (as we may) that $\sum_j\gc_{ij}=1$ for each $i$.
The
proof of Theorem \ref{MainThm}
resembles that of Theorem 33 in \cite{DR}, and
is based on
\begin{obs} \label{L1} For any $i \in [n]$, $k \in \Nn$, and event $Q$ determined by $(\gs^{-1}(j):j \neq i)$,
$$\Pr(B_i=k+1,Q)=\frac{1}{k+1} \sum_{l \in [m]} \gc_{li} \Pr\nolimits^{[m] \sm \{l\}} (B_i=k,Q).$$
\end{obs}
(Recall $\Pr^L$ was defined at the end of Section \ref{Intro}.)
We also use the trivial
\beq{L2}
\min_{i} \frac{\ga_i}{\gb_i} \leq \frac{\ga_1 + \cdots + \ga_k}{\gb_1 + \cdots + \gb_k}
\leq \max_{i} \frac{\ga_i}{\gb_i}
\enq
(for all
$\ga_1 \dots \ga_k,\gb_1 \dots \gb_k\geq 0$ with
$ \gb_1 + \cdots + \gb_k >0$,
where, again, $x/0:=\infty$ when $x>0$).

\mn
{\em Proof of Theorem} \ref{MainThm}.
We proceed by
induction on $m$, omitting the easy base cases with $m=1$.
For (a), it's enough to show that the ratio in question does not
increase when we increase a single
entry---w.l.o.g. the $(n-1)$st---of
one of $a,b$.
Thus, by (\ref{L2}), it suffices to show that
$$
\mbox{$\displaystyle
\frac{\Pr(B_n=k+1|B_{n-1}=t,R)}{\Pr(B_n=k|B_{n-1}=t,R)}~$
{\em is nonincreasing in $t$},}
$$
where $R=\{a_j \leq B_j \leq b_j ~\forall j\in [n-2]\}$;
and
by Proposition \ref{propCVX}, this will follow if we show
\beq{I1}
\frac{\Pr(B_n=k+1|B_{n-1}=t+1,R)}{\Pr(B_n=k|B_{n-1}=t+1,R)}
\leq
\frac{\Pr(B_n=k+1|B_{n-1}=t,R)}{\Pr(B_n=k|B_{n-1}=t,R)}
\enq
for all $t$ for which the probabilities appearing in (\ref{I1}) are positive.
(This is the easy reduction of (a) to (\ref{mukl'})
mentioned earlier.)

By Observation \ref{L1} we may write the left side of (\ref{I1}) as
$$
\frac{\sum_{l \in [m]} \gc_{l,n-1}
\Pr\nolimits^{[m] \sm \{l\}}(B_n=k+1,B_{n-1}=t,R)}
{\sum_{l \in [m]} \gc_{l,n-1}
\Pr\nolimits^{[m] \sm \{l\}}(B_n=k,B_{n-1}=t,R)},
$$
which, by (\ref{L2}), is at most
$$
\max_{l \in [m]} \frac{\Pr\nolimits^{[m] \sm \{l\}}(B_n=k+1|B_{n-1}=t,R)}{\Pr\nolimits^{[m] \sm \{l\}}(B_n=k|B_{n-1}=t,R)}.
$$
Thus, setting $Q=\{B_{n-1}=t\}\wedge R$ and
assuming (w.l.o.g.)
that the maximum occurs at $l=m$, we will have (\ref{I1}) if we show
\beq{I2}
\frac{\Pr(B_n=k+1|Q)}{\Pr(B_n=k|Q)}
\geq \frac{\Pr\nolimits^{[m-1]}(B_n=k+1|Q)}
{\Pr\nolimits^{[m-1]}(B_n=k|Q)}.
\enq
Now
$$\frac{\Pr(B_n=k+1|Q)}{\Pr(B_n=k|Q)} =
\frac{\sum_{j \in [n]} \Pr(\gs(m)=j|Q)
\Pr(B_n=k+1|Q,\gs(m)=j)}
{\sum_{j \in [n]} \Pr(\gs(m)=j|Q) \Pr(B_n=k|Q,\gs(m)=j)},$$
so that (\ref{I2}) will follow (again using (\ref{L2})) from
\beq{I3}
\frac{\Pr(B_n=k+1|Q,\gs(m)=j)}{\Pr(B_n=k|Q,\gs(m)=j)} \geq \frac{\Pr\nolimits^{[m-1]}(B_n=k+1|Q)}{\Pr\nolimits^{[m-1]}(B_n=k|Q)}
~~~~\textrm{ for all } j \in [n]
\enq
(where, again, ``for all $j\in [n]$'' really includes only those for which
\mbox{$\Pr(Q,\gs(m)=j)>0$}).

There are three cases to consider.
If $j=n$, the left side of (\ref{I3}) is
$$\frac{\Pr\nolimits^{[m-1]}(B_n=k|Q)}
{\Pr\nolimits^{[m-1]}(B_n=k-1|Q)},$$
which is at least
the right side of (\ref{I3}) by (part (b) of) our
induction hypothesis.
If $j=n-1$, the left side of (\ref{I3}) is
$$\frac{\Pr\nolimits^{[m-1]}(B_n=k+1|B_{n-1}=t-1,R)}{\Pr\nolimits^{[m-1]}(B_n=k|B_{n-1}=t-1,R)},$$
which is at least the right side of (\ref{I3}) by (part (a) of) the induction hypothesis.
Finally, if $j \neq n-1,n$, the left side of (\ref{I3}) is
$$\frac{\Pr\nolimits^{[m-1]}(B_n=k+1|B_{n-1}=t,R^*)}{\Pr\nolimits^{[m-1]}(B_n=k|B_{n-1}=t,R^*)},$$
where $R^*$ is obtained from $R$ by replacing the condition $a_j \leq B_j \leq b_j$ by
the condition $a_j-1 \leq B_j \leq b_j-1$;
again this is at least the right side of (\ref{I3})
by part (a) of the induction hypothesis.

\medskip
We now turn to (b) and set $Q=\{a_j \leq B_j \leq b_j ~~\forall j\in [n-1]\}$.
Then we have, again using Observation \ref{L1} and (\ref{L2}),
\begin{eqnarray*}
\frac{\Pr(B_n=k+1|Q)}{\Pr(B_n=k|Q)} & = & \frac{k \sum_{l \in [m]} \gc_{ln} \Pr\nolimits^{[m] \sm \{l\}}(B_n=k,Q)}
 {(k+1)\sum_{l \in [m]} \gc_{ln} \Pr\nolimits^{[m] \sm \{l\}}(B_n=k-1,Q)} \\
 & \leq & \max_{l \in [m]} \frac{k \Pr\nolimits^{[m] \sm \{l\}}(B_n=k,Q)}{(k+1)\Pr\nolimits^{[m] \sm \{l\}}(B_n=k-1,Q)} \\
 & \stackrel{\rm w.l.o.g.}{=}& \frac{k \Pr\nolimits^{[m-1]}(B_n=k,Q)}
 {(k+1)\Pr\nolimits^{[m-1]}(B_n=k-1,Q)}
\end{eqnarray*}
(noting that we may assume, by Proposition \ref{propCVX}, that $\Pr(B_n=r|Q)$ is
positive for $r \in \{k-1,k,k+1\}$);
so we will be done if we can show
\beq{showb}
 \frac{\Pr\nolimits^{[m-1]}(B_n=k|Q)}
 {\Pr\nolimits^{[m-1]}(B_n=k-1|Q)}
\leq
 \frac{\Pr(B_n=k|Q)}
 {\Pr(B_n=k-1|Q)}.
\enq

Proceeding as in the proof of part (a), we may rewrite
$$
 \frac{\Pr(B_n=k|Q)} {\Pr(B_n=k-1|Q)}
=
\frac{\sum_{j \in [n]} \Pr(\gs(m)=j|Q) \Pr(B_n=k|Q,\gs(m)=j)}
{\sum_{j \in [n]} \Pr(\gs(m)=j|Q) \Pr(B_n=k-1|Q,\gs(m)=j)};$$
so for (\ref{showb}) it is enough to show
that, for each $j\in [n]$,
$$
\frac{\Pr(B_n=k|Q,\gs(m)=j)}
{\Pr(B_n=k-1|Q,\gs(m)=j)} \geq
\frac{\Pr\nolimits^{[m-1]}(B_n=k|Q)}
{\Pr\nolimits^{[m-1]}(B_n=k-1|Q)},
$$
which, as did (\ref{I3}), follows easily from our induction hypothesis
(here we only need to consider the two cases $j=n$ and $j \neq n$).
\qed

\section{A graphical approach}\label{GA}

\medskip
We begin here with a natural and seemingly new
graph theoretic statement
which we regard as the main point of this section.
Given a multigraph $G$ on vertex set $V$ and
$a,b\in \Nn^V$, let
$\Oh(a,b)=\Oh_G(a,b)$ be the set of orientations of $G$ for which
$$\mbox{$d^+(x)\geq a_x~$ and $~ d^-(x)\geq b_x ~~~\textrm{ for all } x\in V$}$$
and $N(a,b)=N_G(a,b)=|\Oh(a,b)|$.
Here $d^+$ and $d^-$ are, as usual, out- and in-degrees.
We will also use $d_x$ for the degree of $x$ in $G$.
Note we regard a loop (at $x$, say)
as having two orientations, each of which
contributes 1 to each of $d^+(x)$ and $d^-(x)$.
\begin{lemma}\label{Glemma}
If $a,b,r,s\in \Nn^V$ satisfy
$$
\mbox{$a\geq r,s~$ and $~a+b\geq r+s$}
$$
(where the inequalities are with respect to
the product order on $\Nn^V$),
then
\beq{Nab}
N(a,b) \leq N(r,s).
\enq
\end{lemma}
Of course the idea is that it's harder to satisfy a set of demands
that always requires large out-degrees than one for which
these requirements are mixed.
For the sake of comparison, let us
also mention the specialization of Corollary \ref{corNLC}
to the present situation:
\begin{cor}\label{gphcor}
If $a+b= r+s$ then
\beq{Nabcd}
N(a,b)N(r,s)\geq
N(a\vee r,b\wedge s)N(a\wedge r,b\vee s).
\enq
\end{cor}
{\em Proof.}
Interpret vertices of $G$ as urns and edges as
balls, and assume that for each edge (ball) $e$ we have
$\gc_{ex}=1$ or $0$ according to whether $x$ is or is not
an end of $e$.  Then (\ref{Nabcd}) is just
Corollary \ref{corNLC} with $f=d-a-b$ ($= d-r-s$),
where $d = (d_x:x\in V(G))$ is the vector of degrees.\qed

\mn{\em Remark.}
It's possible to simplify the proof of
Lemma \ref{Glemma} using Corollary
\ref{gphcor}; but of course this
depends on Theorem \ref{MainThm},
so is really harder
than the following direct proof.
On the other hand, it's not too hard to derive
Theorem \ref{MainThm}(a) from Lemma \ref{Glemma};
see \cite{Mike}.
(And below we use 
Lemma \ref{Glemma} to prove (\ref{mukl'''}),
which is stronger than Theorem \ref{MainThm}(b).)


\mn
{\em Proof of Lemma \ref {Glemma}.}
We proceed by induction on
$\varphi(G,a,b):=|E(G)|+\sum _{x\in V}(d_x-a_x-b_x)$,
calling
$x \in V$ {\em saturated} if $a_x+b_x=d_x$.
Since $N$ is nonincreasing in each of its arguments,
we may assume $a+b=r+s$ (or we can increase $r$ or $s$).

Suppose first that there
is at least one saturated vertex, $x$.
We may assume there are no loops at $x$,
since otherwise (\ref{Nab})
follows easily from the induction hypothesis
applied to the graph gotten from $G$ by deleting
such loops.
Let $\ga=a_x,\gb=b_x,\rho=r_x,\gs=s_x$,
and let $X=\{e_1 \dots e_{\ga+\gb}\}$
be the set of edges incident with $x$.

Consider a set $\pi$ consisting of $\gb$ pairs $\{e_i,e_j\}\sub X$,
with the $2\gb$ edges appearing in $\pi$ distinct, and, say,
$y_i$ the vertex joined to $x$ by $e_i$
(so the $y_i$'s need not be distinct).
Let $G({\pi})$ be the graph with
vertex set $V \sm \{x\}$ and
edge set $E(G) \sm X \cup \{e_{ij}:\{e_i,e_j\} \in \pi\}$,
where $e_{ij}$ joins $y_i$ and $y_j$.
Let $U(\pi)$ be the set of edges in $X$
not belonging to pairs from $\pi$,
and $U_z(\pi)$ the set of edges of $U(\pi)$ incident to $z$.

Define $a^{\pi},b^{\pi} \in \Nn^{V(G(\pi))}$ by
$$a^{\pi}_z=a_z \textrm{ and } b^{\pi}_z=
\max\{b_z-|U_z(\pi)|,0\} \textrm{ for all } z \in V \sm \{x\} ~~(=V(G(\pi))).$$
For each $\pi$ as above and
$T\in \C{U(\pi)}{\rho-\gb}$
(where
$\C{A}{k}=\{B \sub A:|B|=k\}$),
define $r^{\pi,T},s^{\pi,T} \in \Nn^{V(G(\pi))}$ by
$$r^{\pi,T}_z=\max\{r_z-|U_z(\pi) \sm T|,0\}$$
and
$$s^{\pi,T}_z=\max\{s_z-|U_z(\pi) \cap T|,0\}$$
for all $z \in V \sm \{x\}$.

Each $\gs\in \Oh_{G(\pi)}(a^{\pi},b^{\pi})$ maps naturally
to a (unique) $\hat{\gs}\in \Oh_G(a,b)$, namely:
$\hat{\gs}$ agrees with $\gs$ on $E(G-x)$;
orients all edges of $U(\pi)$ away from $x$; and
orients $e_i$ from $y_i$ to $x$ and $e_j$ from $x$ to $y_j$
whenever $\gs$ orients $e_{ij}$ from $y_i$ to $y_j$
(where, when $y_i=y_j$, we interpret one orientation of the
loop $e_{ij}$ as $y_i\ra y_j$ and the other as $y_j\ra y_i$).
Since each $\tau \in \Oh_G(a,b)$
is in the range of this map for exactly $\C{\ga}{\gb}\gb!$
choice of $\pi$, we have
\beq{NabDec}
N(a,b)=\frac{1}{\C{\ga}{\gb} \gb!} \sum_{\pi} N_{G(\pi)}(a^{\pi},b^{\pi}).
\enq
Similarly,
\beq{NrsDec}
N(r,s)=\frac{1}{\C{\rho}{\gb} \C{\gs}{\gb} \gb!}
\sum_{\pi} \sum_{T \in \C{U(\pi)}{\rho-\gb}} N_{G(\pi)}(r^{\pi,T},s^{\pi,T}).
\enq
Since
$a^{\pi}+b^{\pi} \geq r^{\pi,T}+s^{\pi,T}$
and $a^{\pi} \geq r^{\pi,T},s^{\pi,T}$,
it follows from the induction hypothesis that
\beq{NGpi}
\mbox{$N_{G(\pi)}(a^{\pi},b^{\pi}) \leq N_{G(\pi)}(r^{\pi,T},s^{\pi,T})$
for all $\pi$ and
$T \in \C{U(\pi)}{\rho-\gb}$.}
\enq
(Note
that
$\varphi(G(\pi),a^{\pi},b^{\pi})<\varphi(G,a,b)$,
since $|E(G(\pi))|<|E(G)|$
and, for $z \in V \sm \{x\}$,
$a^{\pi}_z+b^{\pi}_z \geq a_z+b_z-|U_z(\pi)|$, while
the degree of $z$ in $G(\pi)$ is
$d_z-|U_z(\pi)|$.)
Combining (\ref{NabDec}), (\ref{NrsDec}) and (\ref{NGpi}), we have
$$N(r,s) \geq \frac{1}{\C{\rho}{\gb} \C{\gs}{\gb} \gb!}
\sum_{\pi} \C{\ga-\gb}{\rho-\gb} N_{G(\pi)}(a^{\pi},b^{\pi})
= \frac{\ga!\gb! }{\rho!\gs!}
N(a,b)\geq N(a,b),$$
where the last inequality follows from the assumptions
$\ga \geq \rho,\gs$ and $\ga+\gb=\rho+\gs$.

So we may assume there are no saturated vertices.
In this case we fix $x \in V$ with $a_x>b_x$.
(Of course if
there is no such vertex, then $a=b=r=s$ and (\ref{Nab}) is an equality.)
For $\gc,\gd\in \Nn$ let $N'(\gc,\gd)$ be the number of orientations
of $G$ with
$$
(d^+_y,d^-_y)\geq \left\{\begin{array}{ll}
(a_y,b_y)  &\mbox{if $y\neq x$}\\
(\gc,\gd)  &\mbox{if $y=x$,}
\end{array}\right.
$$
and let $N''(\gc,\gd)$ be defined analogously with $(r,s)$
in place of $(a,b)$.
Let
$\ga=a_x,\gb=b_x,\rho=r_x,\gs=s_x$,
so that (\ref{Nab})
is
\beq{NabTS}
N'(\ga,\gb)\leq N''(\rho,\gs).
\enq
By induction we have
\beq{ind}
\mbox{$N'(\gc,\gd)\leq N''(\eta,\xi)$ whenever
$\gc\geq \eta,\xi$ and $\gc+\gd =\eta+\xi > \ga+\gb$.}
\enq
%
We apply this to the identity
\beq{NPgagb}
N'(\ga,\gb)=N'(\ga,\gb+1)+N'(d_x-\gb,\gb).
\enq
If $\ga>\gs$, then, by (\ref{ind}),
the right side of (\ref{NPgagb}) is at most
$$N''(\rho,\gs+1)+N''(d_x-\gs,\gs)=N''(\rho,\gs).$$
If $\ga>\rho$, then, again
using (\ref{ind}),
the right side of (\ref{NPgagb}) is at most
$$N''(\rho+1,\gs)+N''(\rho,d_x-\rho)=N''(\rho,\gs).$$
(And, since $\ga>\gb$, we have at least one of $\ga>\gs$, $\ga>\rho$.)
\qed

\bigskip
The next result isolates (and generalizes)
the main point
in the derivation of
(\ref{mukl'''})
from Lemma \ref{Glemma}.
%
We consider a hypergraph  $\h=\h_1\cup \h_2$ on a set $W$ of size $2l$,
where

\mn
(i)  the edges of $\h_1$ are pairwise disjoint
and

\mn
(ii)
the
edges of $\h_2$ are of size 2 and pairwise disjoint.

\mn
Let $S$ be the set of vertices of $\h$ not covered by
edges of $\h_2$, and $|S|=2t$.
Given $\ga:\h_1\ra \Nn$,
let $N_i$ be the number of partitions $(X,Y)$ of $W$ with each of
$X,Y$ a vertex cover of $\h_2$, each of $|X\cap H|$, $|Y\cap H|$
at least $\ga_{_H}$ for each $H\in \h_1$,
and $|X|=i$.

\begin{lemma}\label{hyplemma}
In the above situation, $N_l\geq \frac{t+1}{t}N_{l+1}$.
\end{lemma}
{\em Proof.}
For $i\in \Nn$ and $\pi$ a
collection of $t-1$ disjoint
2-sets contained in $S$,
let
$\N_i(\pi)$ be the set of partitions as above for
which each of $X,Y$ also covers the edges of $\pi$,
and set $N_i=|\N_i|$.
We assert that (for each $\pi$)
\beq{Npi}
N_l(\pi)\geq 2N_{l+1}(\pi).
\enq
This implies the proposition since (as is easily seen)
$$N_l = \frac{1}{t\cdot t!} \sum_{\pi}N_l(\pi)
$$
and
$$N_{l+1} = \frac{1}{\C{t+1}{2} (t-1)!} \sum_{\pi}N_{l+1}(\pi).
$$

For the proof of (\ref{Npi})
let $x,y$ be the two vertices of $W$ not contained in members
of $\h_2':=\h_2\cup \pi$.
\sugg{
, and set $V'=V\sm \{x,y\}$.
For $X, Y$ as above,
we may identify the partition $V'=(X\cap V')\cup (Y\cap V')$
with the orientation of
$\h_2'$ that directs $\{u,v\}$ from $u $ to $v$
if $u\in X$ (so $v\in Y$).
}
Noting that $(X,Y)\in \N_{l+1}(\pi)$ implies $\{x,y\}\in X$,
we may regard $(X,Y)$ as an orientation of $\h_2'$,
where orienting $\{u,v\}$ from $u$ to $v$ corresponds to putting
$u$ in $X$ (and $v$ in $Y$).
The orientations corresponding to $(X,Y)$'s from $\N_{l+1}$ are those
for which, for each $H\in \h_1$,
$$
\mbox{$d^+(H)\geq \ga_{_H}-|H\cap\{x,y\}|~~$ and $~~d^-(H)\geq \ga_{_H}$,}
$$
where, for the given orientation, $d^+(H)$ (resp. $d^-(H)$) is the number
of oriented edges whose tails (resp. heads) lie in $H$.

If we let $G$
be the multigraph gotten from $(\h,\pi)$
by collapsing each $H\in \h_1$ to a single vertex
(so for example, any $\{u,v\}\in \h_2'$ contained in some $H\in \h_1$
becomes a loop in $G$), then the above discussion says that
$N_{l+1}(\pi) = N_G(a,b)$ (see Lemma \ref{Glemma} for the notation),
where
$$
\mbox{$a_z=
\ga_{_H}-|H\cap\{x,y\}|~~$ and $~~b_z= \ga_{_H}$}
$$
if $z$
is the vertex of $G$
corresponding to $H\in \h_1$, and $a_z=b_z=0$ if $z$ is not of this type
(i.e. $z\in W\sm \cup\{H:H\in \h_1\}$).

A similar discussion shows that $N_l(\pi)=N_G(r,s)+N_G(s,r)$,
where
$$
\mbox{$r_z=
\ga_{_H}-\1_{\{x\in H\}}~~$ and $~~s_z= \ga_{_H}-\1_{\{y\in H\}}$}
$$
if $z$
is the vertex of $G$
corresponding to $H\in \h_1$, and $r_z=s_z=0$ if $z$ is not of this type.
(For example, $N_G(r,s)$ counts pairs $(X,Y)$ with $x\in X$
(and $y\in Y$).)

Finally, Lemma \ref{Glemma} gives
$N_G(r,s), N_G(s,r)\geq N_G(a,b)$, so we have (\ref{Npi}).
(Strictly speaking we may be applying Lemma \ref{Glemma}
with some negative entries in $b,r$ and/or $s$;
but it's easy to see that this slightly more general
version follows from the lemma as stated.)\qed

\medskip
As mentioned earlier, the proof of (\ref{mukl'''})
also requires Theorem \ref{CAPPULC} below.
(If we just wanted (\ref{mukl''}) then
Lemma \ref{hyplemma} alone would suffice.)
For
$\mu\in \m_m:=\m(2^{[m]})$, set
\beq{APseq}
\ga_i(\mu) =
\C{m}{i}^{-1}
\sum\{\mu(A)\mu(\bar{A}):A\sub [m], |A|=i\}
\enq
(where $\bar{A}=[m]\sm A$).
Say $\mu\in \m_{2k}$ has the {\em antipodal pairs property}
(APP) if
\glossary{name=APP,description={A measure $\mu \in \m_{2k}$ has the
{\em antipodal pairs property} if
$$\C{2k}{k}^{-1} \sum_{\eta \in \Omega_{2k},\vert\eta\vert=k}
\mu(\eta)\mu(\uone-\eta)
\geq \C{2k}{k-1}^{-1} \sum_{\eta \in \Omega_{2k},\vert\eta\vert=k-1}
\mu(\eta)\mu(\uone-\eta).$$}}
\mbox{$\ga_k(\mu)\geq \ga_{k-1}(\mu)$},
and $\mu\in \m_m$ has the
{\em conditional antipodal pairs property} (CAPP)
if every measure
obtained from $\mu$ by conditioning on
the values of some $m-2k$ variables (for some $k$) has the APP
(where we view conditioning on the values indexed by $T$
as producing a measure in $\m(2^{[m]\sm T}$).
\begin{thm}[\cite{CAPP}]\label{CAPPULC}
A measure with the CAPP and no internal zeros in
its rank sequence is ULC.
\end{thm}

\mn
{\em Proof of} (\ref{mukl'''}).
It's again enough to show this when $K=[n-1]$, $I=\{n\}$ and $J=\0$.
Setting $U=\gs^{-1}(K)$ and letting $\mu$ be the law of $U$, we
prove the equivalent
\beq{Mu}
\mbox{$\mu$ is ULC.}
\enq

For $A\sub [m]$ and $\gs\in K^A$, say $\gs\in Q$
if it satisfies the conditions in (\ref{Q}), which we now
rewrite
\beq{Q'}
S_j\leq |\gs^{-1}(j)| \leq T_j~~~\forall j \in K,
\enq
where $S_j=a_{t_j}(j)$ and $T_j =a_{t_j+1}(j)-1$.
(This extends the $Q$ of (\ref{Q}), which was
a subset of $[n]^{[m]}$.)
Write $\gs\sim A$ if $\gs\in  K^A\cap Q$ and $\gs\sim l$
if $\gs\sim A$ for some $A$ of size $l$,
and set $T(A) = \prod \{\gc_{in}:i\in [m]\sm A\}$.
Then $\mu$ is given by
$$
\mu(A) ~\propto ~T(A)\sum \{ W(\gs):\gs\sim A\}
~~~~(A\sub [m]).$$

By Theorem \ref{CAPPULC} we will have (\ref{Mu}) if we
show that $\mu$ satisfies the CAPP and its rank sequence
has no internal zeros.
The latter condition is given by Proposition \ref{propCVX},
applied with
$$
f_j = \left\{\begin{array}{ll}
T_j-S_j&\mbox{if $j\in K$}\\
0 &\mbox{if $j=n$}
\end{array}\right.
$$
(so that $\mu(|U|=k)=M_f(S_1 \dots S_{n-1},m-k)$.

To show that $\mu$ has the CAPP,
we should verify the APP for the conditional measures
$\mu_{X,Y} \in \m_{Y \sm X}$ given by
$$\mu_{X,Y}(A) \propto \mu(A \cup X) ~~~~ (A \sub Y \sm X),$$
where
$X \sub Y\sub [m]$ and $|Y \sm X|=2k$ (for some $k$).
Fixing $X,Y$ and letting $A, B$ run over subsets of $Z:=Y\sm X$,
this amounts to
$$
T(X)T(Y) \sum_{|A|=k}\sum\{W(\gs)W(\tau):\gs\sim A\cup X,
\tau\sim (Z\sm A)\cup X\}
~~~~~~~~~~~~~~~~~~~~~~~~~~~~~~~~~
$$
\beq{app1}
~~~~~~~~~~~~~~~~~~~~~~~~~~\geq
\frac{k+1}{k}~T(X)T(Y) \sum_{|B|=k-1}\sum\{W(\ga)W(\gb):\ga\sim B\cup X, \gb\sim (Z\sm B)\cup X\}
\enq
Regard each of $\gs,\tau,\ga, \gb$ in
(\ref{app1}) as a bipartite graph on the vertex set
$[m]\cup K$
 in the natural way.
Then for each pair $(\gs,\tau)$ appearing in (\ref{app1})
the {\em multiset} union
$G=\gs\cup \tau$ is a bipartite multigraph with exactly
$2(|X|+k)$ edges and
$$d_G(i)=\left\{\begin{array}{ll}
2&\mbox{if $i\in X$}\\
1&\mbox{if $i\in Y\sm X$}\\
0&\mbox{if $i\in [m]\sm Y$}
\end{array}\right.
$$
(and similarly for pairs $(\ga,\gb))$.
We may thus rewrite (\ref{app1})
as
$$
\sum_G \sum\{ W(\gs)W(\tau):
\gs\cup \tau=G, \gs\sim |X|+k\}
~~~~~~~~~~~~~~~~~~~~~~~~~~~~~~~~~~~~
$$
$$
~~~~~~~~~~~~~~~~~~~~~~~~
\geq
\frac{k+1}{k}
\sum_G \sum\{ W(\ga)W(\gb)
:\ga\cup \gb=G, \ga\sim |X|+k-1\},
$$
and it is enough to show that for each fixed $G$
we have the corresponding inequality for the inner
sum, i.e.
\beq{stab'}
\sum\{ W(\gs)W(\tau)
:\gs\cup \tau=G, \gs\sim |X|+k\}
\geq
\frac{k+1}{k}\sum\{ W(\ga)W(\gb)
:\ga\cup \gb=G, \ga\sim |X|+k-1\}.
\enq
This has the advantage that the weights no longer play
a role, since for $\gs\dots\gb$ as in
(\ref{stab'}),
$$
W(\gs)W(\tau) = W(\ga)W(\gb);
$$
so we will have (\ref{stab'}) if we show
\beq{NAB}
N_k\geq \mbox{$\frac{k+1}{k}$}N_{k-1},
\enq
where $N_i= N_i(G)$ is the number of
partitions $G=\gc\cup \gd$ with $\gc\sim |X|+i$.

Now let $\h=\h_1\cup \h_2$ be the hypergraph with vertex set
$W=E(G)$,
$$
\h_1=\{H_j:j\in K\},
$$
where $H_j= \{e\in W:j\in e\}$, 
and
$$
\h_2 = \{\{e,f\}:\mbox{$e\neq f$, $e$ and $f$ have the same end in $[m]$}\}.
$$
Define $\ga:\h_1\ra\Nn$ by
$$\ga_{_{H_j}} = \max\{S_j,d_G(j)-T_j\}$$
(recall $S_j,T_j$ were defined following (\ref{Q'})).
We are then in the situation of Lemma \ref{hyplemma}:
a partition $W=W_1\cup W_2$ with $|W_1|=|X|+i$ (and $|W_2|=|Y|-i$)
as in the lemma
(i.e. with $(W_1,W_2)$ in place of $(X,Y)$)
is
the same thing as a partition $G=\gc\cup \gd$ with
$\gc\sim |X|+i$ (and $\gd\sim |Y|-i$),
and the $t$ in the lemma is equal to the present $k$;
so we have (\ref{NAB}).\qed

\nin
{\em Remark.}
Log-concavity results being of some interest, we mention one
appealing specialization of (\ref{mukl'''}).
(See e.g.
\cite{Stanley} or \cite{Brenti} for much more on log-concavity
in combinatorial settings.)
For a bipartite graph $G=(V\cup K,E)$, define a {\em G-map}
to be a function $f:A\ra K$ with $A\sub V$ and
$(v,f(v))\in E \textrm{ for all } v\in A$.
Given $l,u \in \Nn^K$ (with $l_j \leq u_j$),
call a $G$-map {\em valid} if
$|f^{-1}(j)| \in [l_j,u_j]$
for all $j \in K$.
Let $s_k=s_k(G,l,u)$ be the number of
valid $G$-maps $f:A\ra K$ with $|A|=k$.
\begin{thm}\label{TGlu}
For any $G,l,u$ the sequence $(s_0 \dots s_{|V|})$ is ultra-log-concave.
\end{thm}
In the special case $l \equiv 0$, $u \equiv1$, $s_k$ becomes
$\Phi_k=\Phi_k(G)$,
the number of matchings of size $k$ in $G$.
Heilman and Lieb \cite{HL1,HL2} and Kunz \cite{Kunz}
(see also \cite[Chapter 8]{LovPlum})
proved that for any (not necessarily bipartite) graph $G$,
the {\em matching generating polynomial}
$$p(x)=\sum_{k=0}^{\nu} \Phi_k x^k$$
(where, as usual,
$\nu=\max\{k:\Phi_k > 0\}$ is the matching number of $G$)
has all real (negative) roots.
This implies, by Newton's inequalities (e.g. \cite[Theorem 51]{HLP}),
that
\beq{HL}
\mbox{$(\Phi_0 \dots \Phi_{\nu})$ is ULC,}
\enq
which, if $\nu<|V|$,
is somewhat stronger than this case of Theorem \ref{TGlu}.
In contrast, for general $l$ and $u$ the polynomial
$$\sum_{k=0}^{|V|} s_k x^k$$
need not have all real roots.
(For example, let $V=\{y,v,w\}$, $K=\{1\}$,
$E(G)=\{\{y,1\},\{v,1\},\{w,1\}\}$, $l_1=1,u_1=3$.)

Actually Theorem \ref{CAPPULC} can be used to show that
for any $G$,
$(\Phi_0 \dots \Phi_{\tau})$ is ULC
(where, as usual, $\tau=\tau(G)$ is the vertex cover number),
which in particular recovers (\ref{HL}) when $G$ is bipartite.
As this doesn't use Lemma \ref{Glemma}, we won't go into it here.
It would be very interesting to see a combinatorial
proof of (\ref{HL}) for general $G$.

\medskip
Before closing this section we point out one further
consequence of Lemma \ref{Glemma}, which seems to us interesting
for its own sake.
With notation as in the above proof of (\ref{mukl'''}), set
$f(A) = \sum\{W(\gs):
\gs\in K^A\cap Q\}$ ($A\sub [m]$).  We assert that
$f$ satisfies the negative lattice condition:
\beq{NLCf}
f(A\cup B)f(A\cap B)\leq f(A)f(B) ~~~~~\forall A,B\sub [m].
\enq
%
While we don't see how to get (\ref{mukl''}) (or (\ref{mukl'''}))
from this in general, it's not hard to see that it does imply
(\ref{mukl''}) in case the $\gs(i)$'s are i.i.d.,
so gives yet another proof Theorem \ref{Uthm'}.
We omit the details.

\bn
{\em Proof of} (\ref{NLCf}).
%
We may rewrite the inequality as
\beq{WAB'}
\sum\sum\{ W(\gs)W(\tau)
:\gs\sim A\cup B,\tau\sim A\cap B\}
\leq
\sum\sum\{ W(\ga)W(\gb)
:\ga\sim A,\gb\sim B\}.
\enq
As before we regard $\gs,\tau,\ga, \gb$ in
(\ref{WAB'}) as bipartite graphs on
$[m]\cup K$.
For each pair $(\gs,\tau)$ appearing in (\ref{WAB'}),
the (multiset) union
$G=\gs\cup \tau$ is a bipartite multigraph
with
$$d_G(i)=\left\{\begin{array}{ll}
 2&\mbox{if $i\in A\cap B$}\\
 1 &\mbox{if $i\in A\triangle B$}\\
 0 &\mbox{otherwise}
 \end{array}\right.
$$
(and similarly for pairs $(\ga,\gb))$,
and it's enough to show that, for each such $G$,
(\ref{WAB'}) still holds if we restrict to pairs
$(\gs,\tau)$ and $(\ga,\gb)$ with
\beq{stab}
\gs\cup \tau = \ga\cup\gb =G.
\enq
Again the weights ($W(\gs)$ etc.)
cancel and it's enough to show
\beq{NAB'}
N(A\cup B,A\cap B)\leq N(A,B),
\enq
where, for $C,D\sub [m]$, $N(C,D)= N_G(C,D)$ is the number of
partitions $E(G)=\gc\cup \gd$ with $\gc\sim C$ and $\gd\sim D$.

Notice now that we are really counting partitions
$\hat{\gs}\cup \hat{\tau}$ and $\hat{\ga}\cup \hat{\gb}$ of the edges of
$G':=G[(A\cap B)\cup K]$,
since for any $\gs\dots \gb$ (as in (\ref{WAB'}))
satisfying (\ref{stab}),
any edge of $G$ with an end in $A\sm B$ (resp. $B\sm A$)
must belong to $\gs\cap \ga$ (resp. $\gs\cap \gb$).

For $j\in K$ and $C\sub [m]$ write $d_C(j)$ for the number of
edges of $G$ joining $j$ to $C$.
In terms of $\hat{\gs}\dots \hat{\gb}$ the requirement that
$\gs\dots \gb$ satisfy (\ref{Q'})
becomes the condition that for each $j\in K$,
%
\beq{rqmts}\begin{array}{cc}
\mbox{$S_j-d_{A\Delta B}(j)\leq |\hat{\gs}^{-1}(j)|\leq T_j-d_{A\Delta B}(j)$;
$~S_j\leq |\hat{\tau}^{-1}(j)|\leq T_j$;}\\
\mbox{$S_j-d_{A\sm B}(j)\leq |\hat{\ga}^{-1}(j)|\leq T_j-d_{A\sm B}(j)~$ and
$~S_j-d_{B\sm A}(j)\leq |\hat{\gb}^{-1}(j)|\leq T_j-d_{B\sm A}(j)$.}
\end{array}
\enq
Now let $H$ be the multigraph on vertex set $E(G')$ with edge set $\{e_x:x\in A\cap B\}$,
where $e_x$ joins the two edges of $G'$ containing $x$.
We may identify a partition $E(G')=\gc\cup \gd$ with the
orientation of $H$ gotten by directing $e_x$ from $a$ to $b$
whenever $a\in \gc$ and $b\in \gd$,
where $a,b$ are the edges on $x$ in $G'$.
The orientations corresponding to pairs $(\hat{\gs},\hat{\tau})$ as in (\ref{rqmts})
are then those satisfying
$$
\mbox{$S_j-d_{A\Delta B}(j)\leq d^+(j)\leq T_j-d_{A\Delta B}(j)$
and
$~S_j\leq d^-(j)\leq T_j$ $~\forall j\in K$}$$
while those corresponding to pairs $(\hat{\ga},\hat{\gb})$ are those with
$$
\mbox{$S_j-d_{A\sm B}(j)\leq d^+(j)\leq T_j-d_{A\sm B}(j)~$ and
$~S_j-d_{B\sm A}(j)\leq d^-(j)\leq T_j-d_{B\sm A}(j)$ $~\forall j\in K$.}
$$
That the number of orientations of the first type is at most
the number of the second type is then an instance of
Lemma \ref{Glemma}.\qed

\section{Final remarks}\label{Remks}

The most interesting question left open by the present work is whether
Theorem \ref{Uthm} (even without thresholds) extends
to nonidentical balls; that is,
\begin{question}\label{QCNA}
Are generalized urn measures (or generalized threshold
or interval urn measures) CNA?
\end{question}
As mentioned in the introduction, Dubhashi and Ranjan \cite{DR}
showed NA for the
$\xi_{ij}$'s defined in (\ref{xiij}), which immediately gives
NA for generalized threshold urn measures.
That the weaker (than CNA) CNC, at least, does hold for generalized
threshold (or, more generally, ``interval") urn measures
is a special case of the following result; this is a somewhat more
general version of Corollary 34 of \cite{DR}, which implies
CNC for generalized threshold urn measures.
We again take $\Pr$ to be the measure on $[n]^{[m]}$ corresponding
to some $\gc$, and,
for $\A\sub 2^{[m]}$ and
$a,b\in \Nn^{n-1}$,
set
$
p(\A,a,b)=\Pr(\gs^{-1}(n)\in \A |B_j\in [a_j,b_j] ~\forall j\in [n-1]).
$
\begin{thm} \label{DR26}
For any increasing $\A$, $p(\A,a,b)$
is decreasing
in $(a,b)$.
\end{thm}
The proof is more or less the same as that of Theorem
\ref{MainThm}(a), so will not be given here;
see \cite[Theorem 1.23]{Mike}.
(The proof of Corollary 34 in \cite{DR}
is not quite correct,
since
it depends on the incorrect Proposition 24.)

Thus one reason to be interested in whether generalized urn measures
are CNA is that
a negative answer
would provide a counterexample
to an
important conjecture
of Pemantle \cite{Pem} stating
that CNC implies CNA.
(He also conjectures that the Rayleigh property NC+ implies NA$+$.)
Pursuing this a little further,
say $\mu\in \m(\{0,1\})^n$ is $R^+$ if
every $W\circ \mu$ with $W_i\in \{0\}\cup [1,\infty) ~\forall i$ is NC.
An easy simulation shows
that CNC for the class of generalized urn measures
is the same as $R^+$ for this class.
(Note this is without thresholds; it's easy to see that
$R^+$ need not hold for (even ordinary) threshold urn measures.)
So failure of CNA here
would in fact disprove
\begin{conj}
$R^+$ implies CNA,
\end{conj}
a weakening of the first conjecture of
Pemantle above.

At this writing we can (e.g.) give a positive answer to Question \ref{QCNA}
when each ball chooses from just two urns;  this is of course quite special,
but seems of some interest since it corresponds to
in- and out-degree statistics for a random orientation of a graph
(where edges are oriented independently, but the two orientations of
an edge may have different probabilities).
Even this special case seems to require an interesting argument, but we will
not give this here as the paper seems long enough without it.

As far as we can see, even the following very general statement
could be true.
\begin{question}\label{QQ}
Suppose $T_0\cup T_1\cup\cdots \cup T_s$ is a partition of $[m]\times [n]$,
and $a_r, b_r\in \Nn$ for $r=1\dots s$.  Is it true
that the $\xi_{ij}$'s in {\em (\ref{xiij})}
are NA given
$$\{\xi(T_r) \in [a_r,b_r]~\forall r\in [s]\}$$
(where
$\xi(T) =\sum_{(i,j)\in T}\xi_{ij}$)?
\end{question}
This would be a considerable strengthening of CNA for
generalized threshold urn measures.

Let us also just mention one possible approach to
Question \ref{QCNA}.
Recall that for $\mu,\nu\in \m(\{0,1\}^n)$,
$\mu$ {\em stochastically dominates} $\nu$
(written $\mu\succeq \nu$) if
$\mu(\A)\geq \nu(\A)$ for each increasing
$\A\sub \{0,1\}^n$,
and that $\mu$ has the {\em normalized matching property}
if, with $X$ chosen according to $\mu$ and $\xi=|X|$,
 $\mu(\cdot |\xi=k)$ is stochastically
increasing in $k$
(meaning, of course, that  $\mu(\cdot |\xi=k)\succeq  \mu(\cdot |\xi=l)$
whenever $k>l$).
It is not too hard to show
(this is somewhat like the derivation of CNA from CNC
in \cite{FM}) that CNA for generalized interval urn measures
would follow from a positive answer
to

\begin{question}
Is it true that for any
$\gs\in [n]^{[m]}$ with law given by {\rm (\ref{tauW})} and
$Q$ as in {\rm (\ref{Q})},
the law of $\gs^{-1}(K)$ given Q
has the normalized matching property?
\end{question}
%

\medskip
Finally we turn to the conjectures of Farr
(unpublished circa 2004; see \cite{Welsh})
and Welsh \cite{Welsh} mentioned at the end of Section \ref{Intro}.
To put
these in our framework, we add an urn $\gL$
and assume
$$\mbox{$\Pr(\gs(i) = j) = p$ $~~~\forall i\in [m] ,j\in [n]$.}$$
(So $\Pr(\gs(i)=\gL) = 1-np$.)
Let $\I\sub 2^{[m]}$ be decreasing
and set $\A_j=\{\gs^{-1}(j)\in \I\}$
and $\A_J=\cap\{\A_j:j\in J\}$.
Then Farr's conjecture (somewhat rephrased) is

\begin{conj}\label{Farr}
If G is a graph on $[m]$ and
$\I$ is the collection of independent sets of $G$,
then for any disjoint $I,J,K\sub [n]$,
$\A_I\downarrow\A_J$ given $\A_K$.
\end{conj}
It's not clear why this should require that $\I$
be of the type described, and Welsh's conjecture was
that the same conclusion
holds for an arbitrary $\I$.
Here we sketch a counterexample to this stronger version.
At present we don't see how to extend to a counterexample
to Conjecture~\ref{Farr}, though we feel that this too is likely to
be false.

\mn
{\em Example.}
Let $n=3$ and $p=1/3$ (so we don't need $\gL$).
Let $M\cup A\cup B\cup C$ be a partition of $V:=[m]$
with $|M|=s$ (large) and $|A|=|B|=|C|=t=s+3$.
Let
$\I_1=2^{V\sm M}$,
$$\I_2=
\{X\sub V:  \mbox{$|X\cap M| <  .4|M|$ and $X$ meets
at most two of $A,B,C$}\},
$$
and $\I\ =\I_1\cup \I_2$.
Then,
we assert,

$$
\Pr(\A_{\{3\}})\Pr(\A_{\{1,2,3\}}) > \Pr(\A_{\{1,3\}})\Pr(\A_{\{2,3\}}),
$$
which contradicts Welsh's conjecture (with $I=\{1\}$, $J=\{2\}$
and $K=\{3\}$).
We omit the precise calculations; roughly,
with $\ga=(2/3)^t$ and $c=(3/2)^{3}$, we have (as $t\ra\infty$)
$$
\Pr(\A_L)\sim \left\{\begin{array}{cl}
(c+3)\ga&\mbox{if $|L|=1$}\\
6\ga^2&\mbox{if $|L|=2$}\\
6\ga^3&\mbox{if $|L|=3$}.
\end{array}\right.
$$

\bn
{\bf Acknowledgment}
Parts of this work were carried out while the authors were
visiting the Isaac Newton Institute and while the first author
was visiting MIT.  The hospitality of both is gratefully acknowledged.

\end{document}